\documentclass[a4paper,11pt]{article}

\usepackage[margin=1in]{geometry}  

\usepackage{moreverb}

\usepackage[dvips,colorlinks,bookmarksopen,bookmarksnumbered,citecolor=red,urlcolor=red]{hyperref}

\usepackage{amsfonts,amsmath,amssymb,amsthm}
\usepackage{graphicx}

\usepackage{algorithm}
\usepackage{algorithmic}

\newcommand{\RR}{\mathbb{R}}

\newtheorem{theorem}{Theorem}[section]

\begin{document}

\title{Graphical representation of separatrices of attraction basins in two and three dimensional dynamical systems}

\author{Roberto Cavoretto, Alessandra De Rossi, Emma Perracchione and Ezio Venturino\footnote{Department of Mathematics \lq\lq G. Peano\rq\rq, University of Torino, via Carlo Alberto 10, I--10123 Torino, Italy. 
}}

\date{}

\maketitle

\begin{abstract}
In this paper we consider the problem of reconstructing separatrices in dynamical systems. In particular, here we aim at partitioning the domain approximating the boundaries of the basins of attraction of different stable equilibria. We start from the 2D case sketched in \cite{cavoretto11} and the approximation scheme presented in \cite{cavoretto11,C-D-P-V}, and then we extend the reconstruction scheme of separatrices in the cases of three dimensional models with two and three stable equilibria. For this purpose we construct computational algorithms and procedures for the detection and the refinement of points located on the separatrix manifolds that partition the phase space. 
The use of the so-called meshfree or meshless methods is used to reconstruct the separatrices.
\end{abstract}





\section{Introduction}

The objective of this research is the formulation of a new algorithm for the graphical reconstruction of unknown surfaces which partition the three dimensional space into disjoint sets. Even if some techniques to prove the existence of invariant sets have already been developed, none of them, except for particular and well-known cases, allows to have a graphical representation of the separatrix manifolds (see  \cite{Dellnitz,Johnson}).

We extend previous preliminary studies \cite{cavoretto11,C-D-P-V}, to obtain an algorithm, written as a \textsc{Matlab} routine, which is able to reconstruct the boundaries of the basins of attraction of different stable equilibria. The numerical tools that are involved are: (i) a bisection-like routine to detect the points lying on the separatrix manifold; (ii) a reduction scheme to select only the most significant of them, thereby reducing redundancy and incrementing spatial uniformity; and (iii) the interpolation reconstruction based on the partition of unity method. Compactly supported radial basis functions (CSRBFs) are used for the local approximations (see, e.g., \cite{Fasshauer07,Fasshauer11,Iske11,Wendland05}), because they are effective and efficient. Indeed they interpolate accurately and stably large numbers of scattered data (see \cite{cavoretto14,cavoretto14a}).

We now briefly illustrate the importance of having such a versatile tool available in applied sciences.
In mathematical applications to real life problems, dynamical systems constitute a powerful modeling tool.
In general, the trajectories in the phase space usually tend to stable equilibrium configurations,
although persistent oscillations or chaotic behavior can certainly arise for particular parameter settings
\cite {Wiggins03}. When these latter regimes do not arise,
it is often the case that the model in consideration presents more than one possible stable equilibria
for the same parameter values. The system outcome is in such case determined by its present state. In
other words, trajectories with different initial conditions will possibly converge toward different
equilibria, depending on the locations of their respective initial conditions. The set of all points
that taken as initial conditions will have trajectories all tending to the same equilibrium is called
the basin of attraction of that equilibrium point. It is then apparent that, to assess the future
system's behavior, it is of paramount importance the accurate determination of these basins. Nearby
initial points could indeed lead to completely different system's outcomes. Thus in real life phenomena, a detailed knowledge of the attraction basins allows to assess  optimal configurations by bringing the system trajectory into a suitable attraction basin and away from the uncertainty of initial conditions close to the separatrices. In such cases a graphical knowledge of the attraction basins is therefore essential to check if the initial condition lies in a \emph{safe} area.

In this paper we describe efficient methods for the reconstruction of these domains. This is achieved by constructing procedures for the determination of the curves (in 2D) and surfaces (in 3D) that are the boundaries of these basins of attraction.

The paper organization follows. Section \ref{pu} describes the basic tools for the approximation scheme, in particular the Partition of Unity method. The following Section deals with the efficient implementation of the algorithm for the separatrix curves and surfaces. Specifically, once we briefly illustrate for convenience of the reader the algorithm used to reconstruct separatrix curves,  we present the technique for obtaining the surfaces separating two and three stable equilibria in systems with three differential equations. We then relate on the numerical results obtained in all these cases. Finally, the concluding Section summarizes the findings of the paper and analyzes some crucial steps of the algorithm which need further investigations.  
\section{Meshless interpolation methods} \label{pu}

Given a set of data, i.e. measurements and locations at which these measurements were obtained, the aim is to find a function ${\cal I}$ which is a \emph{good} fit to the given data. In the following our criterion for a good fit is that the function $ { \cal I }$ must exactly match to given measurements at the corresponding locations.
This kind of problems arises in many scientific disciplines, such as physics, engineering, finance, biomathematics and medicine. 
Typically, in applied sciences, the measurements taken either at different times or from different sensors or viewpoints are not sampled from a regular and uniform grid, but data are irregularly distributed or scattered. Therefore, the approximation process of the rule which exactly matches the given measurements from irregular locations or data is commonly known as scattered data interpolation.
In order to solve such problem we use meshfree or meshless methods since they allow to work  with a large number of points. 
Furthermore, obviously, they are independent from a mesh and thus they are suited for changes in the geometry of the domain.

\subsection{Radial basis function interpolation} 
In this subsection we briefly review the partition of unity approximation based on a localized use of RBF interpolants.

Given a set ${ \cal P}= \{ \boldsymbol{p}_i \in \mathbb{R}^{s},i=1, \ldots, n \}$ of $n$ distinct \textsl{data points} or \textsl{nodes} in a domain $ \Omega \subseteq \mathbb{R}^{s}$, and a corresponding set $ {\cal F}= \{ f_i = f(\boldsymbol{p}_i)  ,i=1, \ldots, n \}$ of \textsl{data values} or \textsl{function values} obtained by sampling some (unknown) function  $f: \Omega \longrightarrow \mathbb{R}$, the standard RBF interpolation problem is to find an interpolant $R:\Omega \longrightarrow \mathbb{R}$ of the form
\begin{equation}
R( \boldsymbol{p})= \sum_{k=1}^{n} \alpha_k \phi (||  \boldsymbol{p} - \boldsymbol{p}_k  ||_2), \quad  \boldsymbol{p} \in \Omega,
\label{rad1}
\end{equation}
where $||\cdot||_2$ is the Euclidean norm, and $ \phi: [0, \infty) \longrightarrow \mathbb{R}$ is a RBF, \cite{Buhmann03}.

 The coefficients $ \{ \alpha_k \}_{k=1}^{n} $ are determined by enforcing the interpolation conditions
\begin{equation}
R( \boldsymbol{p}_i)=f_i, \quad i=1, \ldots, n.
\label{int1}
\end{equation}
Imposing the conditions \eqref{int1} leads to a symmetric linear system of equations
\begin{equation}
\Phi \boldsymbol{\alpha}= \boldsymbol{f},
\label{sys1}
\end{equation}
where $\Phi_{ik} = \phi (||  \boldsymbol{p}_i - \boldsymbol{p}_k  ||_2)$, $i,k=1, \ldots, n$, $\boldsymbol{\alpha}= [\alpha_1, \ldots, \alpha_n]^T$, and $  \boldsymbol{f} =[f_1, \ldots , f_n]^T$. 

In next subsections we will consider $ \phi$  strictly positive definite, since it guarantees the  interpolation problem is well-posed,
i.e. a solution to the problem exists and is unique. This follows from the fact that, using  strictly positive definite functions
the interpolation matrix is positive definite and hence nonsingular.

\subsection{The partition of unity method}

The basic idea of the partition of unity method is to start with a partition of the open and bounded domain
$ \Omega$ into $d$ subdomains $ \Omega_j$, such that $ \Omega  \subseteq \bigcup_{j=1}^{d} \Omega_j$, with some mild overlap among the subdomains \cite{Melenk96}. The subdomains covering the domain $ \Omega$ are usually supposed to be circular 
patches. Associated with the subdomains we choose partition of unity weight functions $W_j$, i.e. a family of compactly supported, nonnegative and continuous functions subordinate to the subdomain $\Omega_j$, such that $ \sum_{j = 1}^d W_j( \boldsymbol{p})=1$ on $\Omega$ and ${\rm supp}(W_j)  \subseteq \Omega_j$.

Then, for each subdomain $ \Omega_j$, a local interpolation problem in the form \eqref{rad1} is considered and the global approximant is given by:
\begin{equation}
{\cal I}( \boldsymbol{p})= \sum_{j=1}^{d} R_j( \boldsymbol{p} ) W_j ( \boldsymbol{p}), \quad \boldsymbol{p} \in \Omega.
\label{intg}
\end{equation}

Note that if the local approximants satisfy the interpolation conditions  \eqref{int1}, then the
global approximant also interpolates at this nodes since:
\begin{align*}
{\cal I}( \boldsymbol{p}_i ) = \sum_{j=1}^{d} R_j( \boldsymbol{p}_i ) W_j ( \boldsymbol{p}_i) = \sum_{j \in I(\boldsymbol{p}_i)} f( \boldsymbol{p}_i ) W_j ( \boldsymbol{p}_i) = f( \boldsymbol{p}_i),
\end{align*}
where:
\begin{displaymath}
I(\boldsymbol{p}_i)= \{ j / \boldsymbol{p}_i \in \Omega_j \}.
\end{displaymath}
According to \cite{Wendland02}, if we assume to have a $ k$-stable partition of unity, then the derivatives of the weight functions satisfy 
\begin{align*}
||D^{ \mu} W_j ||_{L^{ \infty} ( \Omega_j)} \leq \frac{C_{ \mu} }{ \delta_j^{ | \mu|}}, \quad | \mu |	 \leq k, \quad  \forall \mu \in \mathbb{N}^{s},
\end{align*}
where $ \delta_j$ is the diameter of $\Omega_j$ and $C_{ \mu} > 0$ is a constant. 

Now, after defining the space $C_{\nu}^k(\RR^s)$ of all functions $f \in C^k$ whose derivatives of order $|\mu|=k$ satisfy $D^{\mu}f(\boldsymbol{p})= O(||\boldsymbol{p}||_2^{\nu})$ for $||\boldsymbol{p}||_2 \rightarrow 0$, we have the following approximation theorem (see \cite{Wendland02}).

\begin{theorem}
	 Let $\Omega \subseteq  \RR^s$ be open and bounded and  ${\cal P} = \{\boldsymbol{p}_i, i=1,$ $\ldots,n \}\subseteq \Omega$. Let $\phi \in C_{\nu}^k(\RR^s)$ be a strictly positive definite function. Let $\{\Omega_j\}_{j=1}^{d}$ be a regular covering for $(\Omega, {\cal P})$ and let $\{W_j\}_{j=1}^{d}$ be $k$-stable for $\{\Omega_j\}_{j=1}^{d}$. Then the error between $f \in {\cal N}_{\phi}(\Omega)$, where ${\cal N}_{\phi}$ is the native space of $\phi$, and its partition of unity interpolant (\ref{intg}) can be bounded by
\begin{equation}
	|D^{\mu}f(\boldsymbol{p}) - D^{\mu}{\cal I}(\boldsymbol{p})| \leq C h_{{\cal P}, \Omega}^{(k+\nu)/2 - |\mu|} |f|_{{\cal N}_{\phi}(\Omega)}, 
	\nonumber
\end{equation}
$ \forall \boldsymbol{p} \in \Omega$, $ |\mu| \leq k/2$, and 
	$h_{{\cal P}, \Omega}$ being the so-called \textsl{fill distance}, whose definition is given by 
\begin{equation}
h_{{\cal P}, \Omega} = \sup_{\boldsymbol{p} \in \Omega}\min_{\boldsymbol{p}_j\in {\cal P}} ||\boldsymbol{p}-\boldsymbol{p}_j||_2.
\end{equation} 
\end{theorem}  

We require some additional assumptions on the regularity of $\Omega_j$:
\begin{enumerate}
	\item[(i)] for each $\boldsymbol{p} \in \Omega$ the number of subdomains $\Omega_j$ with $\boldsymbol{p} \in \Omega_j$ is bounded by a global constant $K$;
	\item[(ii)] each subdomain $\Omega_j$ satisfies an interior cone condition (see \cite{Wendland05});
	\item[(iii)] the local fill distances $h_{{\cal P}_j, \Omega_j}$ are uniformly bounded by the global fill distance $h_{{\cal P}, \Omega}$, where ${\cal P}_j={\cal P} \cap \Omega_j$.
\end{enumerate}

This local approach enables us to decompose a large problem into many small problems,
and at the same time ensures that the accuracy obtained for the local fits is carried over to the
global fit. In the reconstruction of the attraction basins the use of compactly supported radial basis functions,
 since  it is well known that they guarantee a good compromise between accuracy and stability, is strongly advised.
 Moreover usually, it can be highly advantageous to work with locally supported 
functions since they lead to sparse linear systems. 
Thus in next subsection  some  widely used CSRBFs will be illustrated.

\subsection{Interpolation with CSRBFs} 
In this subsection we consider the most popular families of CSRBFs such as Wendland's, Wu's and Gneiting's functions, \cite{Fasshauer07,Wendland05}. 
Wendland \cite{Wendland05} found a class of RBFs which are smooth, compactly supported, and strictly
positive definite. They consist of a product of a truncated power function and a low degree polynomial.

Here we list few of the most commonly Wendland's functions:
\begin{equation}
\left.
\begin{array}{rcll}
\varphi_1(r) & = & \left( 1-c r\right)_+^{4} \left(4cr+1\right),                & \hspace{0.3cm} {\rm C^2} \\ 
\varphi_2(r) & = & \left( 1-c r\right)_+^{6}\left(35(c r)^2+18c r+3\right),          & \hspace{0.3cm} {\rm C^4} \\ 
\end{array}
\right.
\end{equation}
where $c \in  \mathbb{R}^{+}$ is the so-called shape parameter. The functions  are  non negative for $r \in [0,1/c ]$. 

Another class of compactly supported and strictly positive definite functions is the class of Wu's functions. Similarly to the Wendland's ones they are obtained  by using the truncated power function. In what follows we report some examples of Wu's functions:
\begin{equation}
\left.
\begin{array}{ll}

\psi_1(r)  =  \left( 1-c r)^{5}_{+} (5 (cr)^4+25 (cr)^3 \right. &   \nonumber \\
 \hspace{1.cm} \left.+48 (cr)^2+40c r+8 \right), &   \hspace{0.3cm} {\rm C^2}\nonumber \\
\psi_2(r) =  \left( 1-cr)^{6}_{+} (5 (cr)^5+30 (cr)^4 \right. &  \nonumber \\
 \hspace{1.cm} \left.  +72(cr)^3+82(cr)^2 +36cr +6\right).  & \hspace{0.3cm} {\rm C^4} 
\end{array}
\right.
\end{equation}

We conclude this subsection by illustating some examples of the class of the so-called Gneiting's or oscillatory  functions.
This family  can be obtained starting from the Wendland's one, (see \cite{Fasshauer07} for details). As example, the following functions are strictly positive definite and radial on $\mathbb{R}^2$:
\begin{equation}
\left.
\begin{array}{rcll}
\tau_1(r) &  = & \left( 1-cr\right)_+^{7/2} \left(-{{135}\over{8}}(c r)^2  +{{7}\over{2}}c r+  1 \right),   & \hspace{0.5cm} {\rm C^2}  \\
\tau_2(r) & = & \left( 1-cr\right)_+^{5} \left(-27(cr)^2+5c r +1 \right).              & \hspace{0.5cm} {\rm C^2}  
 \nonumber \\ 
\end{array}
\right.
\end{equation}

\section{Approximation of separatrices}
We remark that our aim is the reconstruction of  separatrix surfaces partitioning a 3D
phase state into two or more regions. Anyway here, for convenience of the
reader, we illustrate the basic numerical tools used to approximate separatrix curves.
The differences from our previus work, \cite{cavoretto11},
in the reconstruction of separatrix curves 
 will be pointed out. They concern essentially the set of initial conditions used to find the separatrix
points and the refinement algorithm. The latter is employed to find a set of nodes 
well distributed in the domain in which the separatrix curve is defined.

The models used in our examples are a competition-like model for the 3D case and a dynamical system 
modeling  infectious diseases for the 2D case.

\subsection{Calculation of separatrix curves} \label{Curves}

To test the robustness of our  algorithm in a 2D dynamical system, when bistability occurs, we consider the following model describing a population affected by a disease, \cite{Hilker}:
\begin{equation} \label{model2d_hilker}
\begin{array}{ll}
\frac{ \displaystyle  dP}{ \displaystyle  d t}=r(1-P)(P-u)P -\alpha I,  & \textrm{} \\
\vspace{.01cm}\\
\frac{ \displaystyle  dI}{ \displaystyle  d t}= [- \alpha -d-ru+(\sigma-1)P-\sigma I]I, & \textrm{}
\end{array}
\end{equation}
where $P$ is the dimensionless
total population that is composed of infecteds $I$ and
susceptibles $P-I$, (see \cite{Hilker} for further investigations). The Allee effect governs the first equation. It is easy to verify that $E_0=(0,0)$, $E_1=(1,0)$ and $E_2=(u,0)$ are equilibria of the system \eqref{model2d_hilker}, while for the study of the endemic steady states  see \cite{Hilker}.
Here we omit the analytical study of the model, already outlined in \cite{Hilker}.

As suggested by \cite{Hilker} we set $r=0.2$, $u =0.1$, $d=0.25$ and $\alpha=  0.1$; furthermore we fix $\sigma=2.5$. With this choice exists exactly one endemic steady state $E_4 \approx (0.6663, 0.2518)$ which is a stable equilibrium points. Moreover the origin is stable, $E_1$ is an unstable equilibrium point and $E_2=(0.1,0)$  is the saddle point partitioning the phase plane domain.

This situation suggests the existence of a curve separating the paths tending to  disease-free equilibrium point from the trajectories 
tending to the endemic steady state. The separatrix curve divides the phase plane into two subregions, called
basins of attraction of each respective equilibrium \cite{Arrowsmith90,Murray02}. Trajectories originating
in each of them tend to the unique equilibrium $E_0$ or $E_4$ which lies within the basin.

At first, to determine the separatrix curve for (\ref{model2d_hilker}), we need to consider a set of points as initial conditions in a square domain $[0,\gamma]^2$, where $\gamma \in \RR^+$(in the following we will fix $\gamma=10$).
Then we take points in pairs within $[0,\gamma]^2$
and we check whether trajectories from these two points converge to different
equilibria. If this the case, we then proceed with a bisection algorithm along the segment
joining these points in order to determine a separatrix
point. Once we find a set of points on the separatrix,
we perform a refinement of this set,
which is then interpolated using a suitable method (see Section \ref{pu}).

More precisely, in the 2D case we start considering $n$ equispaced initial conditions on each edge of the
square $[0,\gamma]^2$  and the bisection algorithm is applied with the following initial conditions: 
$$
(x_i,0) \quad \textrm{and} \quad (x_i,\gamma), \quad i=1, \ldots, n, 
$$
$$
(0,y_i) \quad \textrm{and} \quad (\gamma,y_i), \quad i=1, \ldots, n. 
$$
In \cite{cavoretto11} the set of initial conditions is taken as a grid of points on the 
	square and for each possible couple of points the bisection-like routine is performed. It is evident that this approach is computationally expensive, compared with the one 
	proposed here which consists in taking points only on the boundary of the square. Moreover in 
	view of performing a similar algorithm in the 3D case such a cheaper (in terms of computational complexity) procedure
	becomes essential. 

Performing the bisection algorithm, a certain number of points  is found on the separatrix curve.
The $N$ points found by the bisection algorithm are collected in a matrix 
$A=(a_{j,k})$, $j=1, \ldots, N$, $k=1,2$, and then refined in order to obtain a smaller set of well distributed nodes on the separatrix curve. 
So we define:
$$
M_x=  \max_{ \displaystyle j} (a_{j,1}),  \quad  j=1, \ldots, N,
$$
$$
M_y=  \max_{ \displaystyle j} (a_{j,2}),  \quad  j=1, \ldots, N,
$$
and we divide  $[0,M_x] \times [0,M_y]$ in $L^2$ subintervals. Then we make an average of the separatrix points on each subintervals.

\begin{figure}[ht!]
\begin{center}
  \includegraphics[height=.24\textheight]{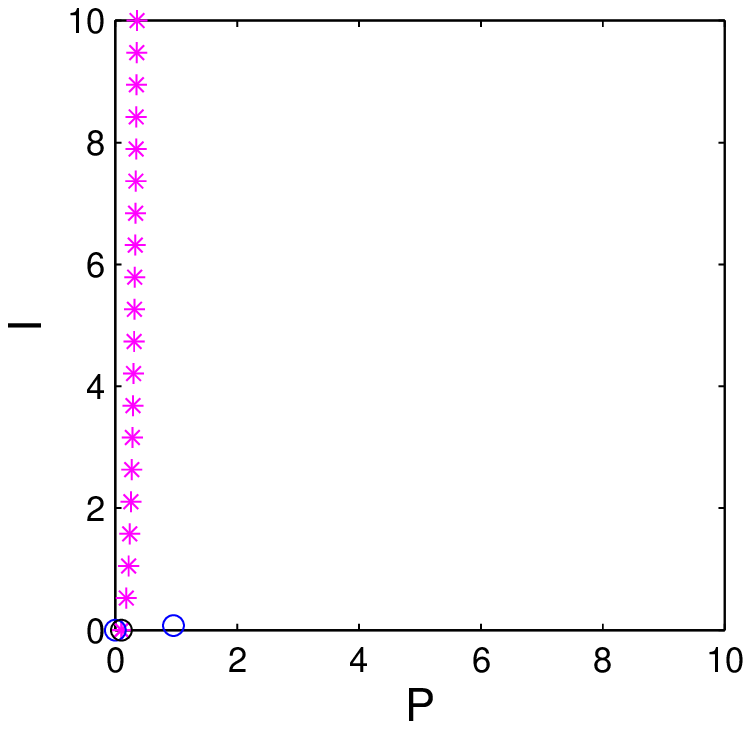} \\
  \includegraphics[height=.24\textheight]{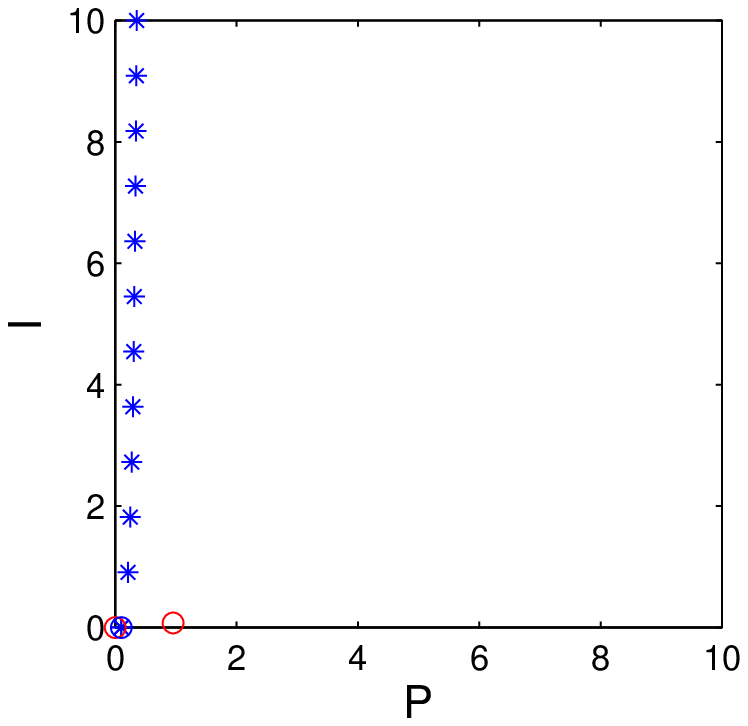} 
  \caption{Set of points detected by the bisection algorithm (top) and
grid of points found by the refinement algorithm (bottom) in the 2D case. Note that these are more
evenly distributed and especially double points or very close points are eliminated.}
\label{figura1}
\end{center}
\end{figure}

For example, Figure \ref{figura1} (top) shows the points found using $n=20$. Choosing $L=12$ and considering the $N=22$ points picked up on the separatrix curve, the refinement process provides us the $K=12$ points reported in Figure \ref{figura1} (bottom). To this set we add the saddle point $E_2$. The refined grid obtained in this way is
then interpolated in order to find the separatrix.

\subsection{Reconstruction of separatrix surfaces} \label{Surf}

Now we consider the following competition model (see \cite{gosso12,gurnell04}):
\begin{equation} \label{model3d}
\begin{array}{ll}
\frac{ \displaystyle  dx}{ \displaystyle  dt}=p \big(1- \frac{ \displaystyle  x}{ \displaystyle  u} \big)x-axy-bxz,  & \textrm{} \\
\vspace{.01cm}\\
\frac{ \displaystyle  dy}{ \displaystyle  dt}=q \big(1- \frac{ \displaystyle  y}{ \displaystyle  v} \big)y-cxy-eyz, & \textrm{} \\
\vspace{.01cm}\\
\frac{ \displaystyle  dz}{ \displaystyle  dt}=r \big(1- \frac{ \displaystyle  z}{ \displaystyle  w} \big)z-fxz-gyz,  & \textrm{} 
\end{array}
\end{equation}
where $x$, $y$ and $z$ denote the three populations, each one competing with both the other ones
in the same environment. We assume that all parameters are nonnegative:
$p$, $q$ and $r$ are the growth rates of the three populations, $a$, $b$, $c$, $e$, $f$ and $g$ denote the competition rates and
$u$, $v$ and $w$  are the carrying capacities of $x$, $y$ and $z$, respectively.

There are eight equilibrium points. The origin $E_0 = (0, 0, 0)$ and the points associated with the survival of only one population $E_1 = (u, 0, 0)$, $E_2 = (0, v, 0)$ and $E_3 = (0, 0, w)$ are always feasible. Then we have the equilibria with two coexisting
populations: \\
\begin{displaymath}
\left.
\begin{array}{ll}
E_4 = & \bigg( \frac{\displaystyle uq(av-p)}{\displaystyle cuva-pq},\frac{\displaystyle pv(cu-q)}{\displaystyle cuva-pq},0\bigg),\\
E_5 = &\bigg( \frac{\displaystyle ur(bw-p)}{\displaystyle fuwb-rp},0, \frac{\displaystyle wp(fu-r)}{\displaystyle fuwb-rp}\bigg),\\
E_6 =& \bigg( 0, \frac{\displaystyle vr(we-q)}{\displaystyle gvwe-qr}, \frac{\displaystyle wq(vg-r)}{\displaystyle gvwe-qr}\bigg).\\
\end{array}
\right.
\end{displaymath}

The feasibility conditions for the point $E_4$ are:
\begin{equation}
\begin{array}{l}
q<cu,  \hspace{0.3 cm}  p<av \quad \textrm{or} \quad q>cu,  \hspace{0.3 cm}  p>av. 
\end{array}
\label{f1}
\end{equation}
Similarly we find the feasibility conditions for $E_5$ and $E_6$ shown in \eqref{e5am1} and \eqref{e6am1}, respectively:
\begin{equation}
\begin{array}{l}
p<bw,  \hspace{0.3 cm}  r<fu \quad \textrm{or} \quad p>bw,  \hspace{0.3 cm}  r>fu, \\
\end{array}
\label{e5am1}
\end{equation}
\begin{equation}
\begin{array}{l}
q<we, \hspace{0.3 cm}  r<vg \quad \textrm{or} \quad  q>we, \hspace{0.3 cm}  r>vg. \\
\end{array}
\label{e6am1}
\end{equation}

Finally we have the coexistence equilibrium:
\begin{displaymath}
\left.
\begin{array}{ll}
E_7 = & \bigg(\frac{\displaystyle u[p(gvwe-qr)-avr(we-q)-bwq(vg-r)]}{\displaystyle p(gvwe-qr)+uva(rc-fwe)+uwb(fq-gcv)},\\
& \frac{ \displaystyle v[q(fuwb-pr)-rcu(wb-p)-pew(fu-r)]}{ \displaystyle q(fuwb-pr)+cuv(ra-gwb)+evw(gp-afu)},\\
&\frac{ \displaystyle r[(cuva-pq)-gpv(cu-q)-ufq(va-p)]}{\displaystyle r(cuva-pq)+bwu(fq-vcg)+evw(gp-fua)} \bigg).
\end{array}
\right.
\end{displaymath}
The stability and feasibility conditions of $E_7$ have been studied only with numerical simulations.

To study the stability conditions for the other equilibria, let us consider the system Jacobian:
\begin{displaymath}
\begin{small}
{\mathcal J_2} =\left[
\begin{array}{ccc}
 \bar{A} & -ax & -bx \\
-cy & \bar{B} & -ey \\
-fz & -gz & \bar{C}\\
\end{array}
\right],
\end{small}
\label{Jac}
\end{displaymath}
where 
\begin{align*}
\begin{array}{rcl}
\bar{A} &=& p \big(1-\frac{ \displaystyle  2x}{ \displaystyle  u} \big)-ay-bz, \\ 
\bar{B} &=& q \big(1-\frac{ \displaystyle  2y}{ \displaystyle  v} \big)-cx-ze, \\ 
\bar{C} &=& r \big(1-\frac{ \displaystyle  2z}{ \displaystyle  w} \big)-fx-gy. 
\end{array}
\end{align*}
Here, for shortness, we omit details and we summarize results in Table \ref{tabella1}. Our studies for the equilibria with two coexisting populations is based on factorizing the characteristic equations to get one eigenvalue and then the  Routh-Hurwitz criterion is applied to the remaining factors of the characteristic equations, i.e.  to  submatrices properly obtained starting from ${\mathcal J_2}$ \cite{gosso12}. All the results, shown in Table \ref{tabella1}, have been verified with symbolic calculations carried out  with \textsc{Maple}.


\begin{table}[!htbp]
\begin{center}
\begin{tabular}{c|c}
\hline
\rule[0mm]{0mm}{3ex}
Eq. & Stability  \\ 
\hline
\rule[0mm]{-1mm}{3ex}
\bfseries $E_0$   &   unstable\\
\bfseries  $E_1$  &  $r<fu$, $q<cu$\\ 
\bfseries  $E_2$   &  $r<vg$, $p<av$\\
\bfseries  $E_3$  &  $q<ew$, $p<bw$ \\
\bfseries  $E_4$  &  $q>cu$, $  p>av$,\\
\bfseries  & $r(cuva-pq)>pvg(cu-q)+ufq(va-p)$ \\
\bfseries  $E_5$  &   $p>bw$,   $r>fu$,\\
\bfseries  & $q(fuwb-pr)>wpe(fu-r)+rcu(wb-p)$ \\
\bfseries  $E_6$  & $ q>we,$ $r>vg$, \\
\bfseries  &  $p(gvwe-rq)>bwq(vg-r)+avr(we-q)$ \\  
\\[\smallskipamount]
\hline
\end{tabular}
\end{center}
\vspace{0.3cm}
\caption{Stability conditions for the equilibria of the system \eqref{model3d}.}
\label{tabella1}
\end{table}
From Table \ref{tabella1} we deduce that for suitable parameters choices the system admits two or three stable equilibria.
For example, with the choice of the parameters $p = 1$, $q =1$, $r = 2$, $a = 1$,
$b=2$, $c =0.3$, $e=1$, $f=3$, $g=2$, $u=1$, $v=0.2$, $w=9.5$, the points $E_3$ and $E_4$ are the only stable
equilibria.  While, with parameters $p = 1$, $q = 2$, $r = 2$, $a = 2$, $b=5$, $c =3$, $e=7$, $f=3$, $g=5$, $u=3$, $v=2$, $w=2$, the points $E_1$, $E_2$ and $E_3$ are  stable
equilibria.
We verify numerically that with these choices, in both cases, $E_7$ is a saddle point. 
This suggests the existence of  separating surfaces partitioning the model domain into two and three basins of
attraction, respectively. The problem of the reconstruction of the surface separating two stable equilibria has been analyzed in \cite{C-D-P-V}, but here we use a different refinement algorithm.
In fact, with this method we are able to reconstruct with one and the same technique
the surfaces separating the basins in the cases of both two and three stable equilibria.

\begin{figure}[ht!]
\begin{center}
 \includegraphics[height=.27\textheight]{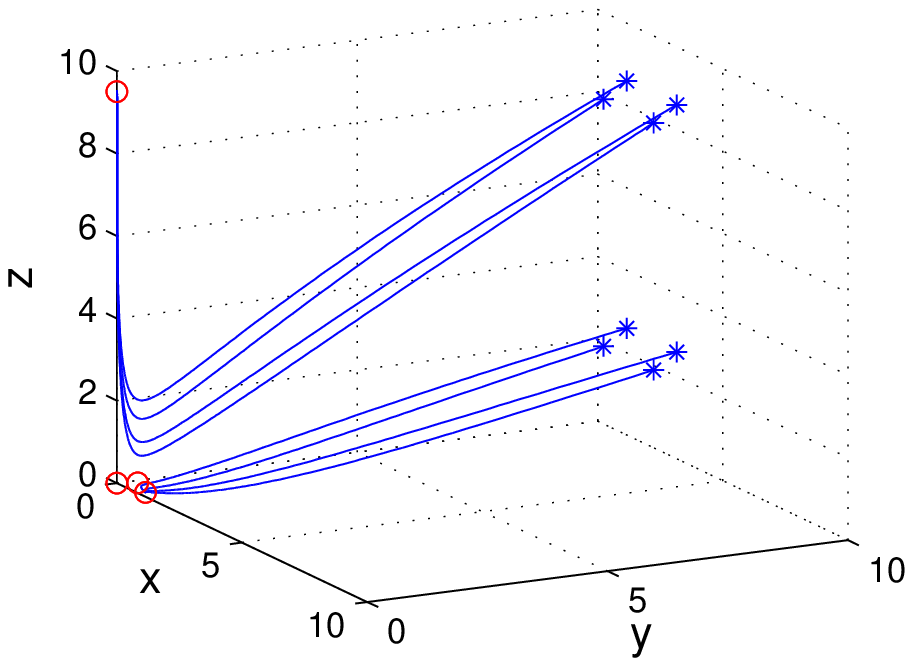}
 \includegraphics[height=.27\textheight]{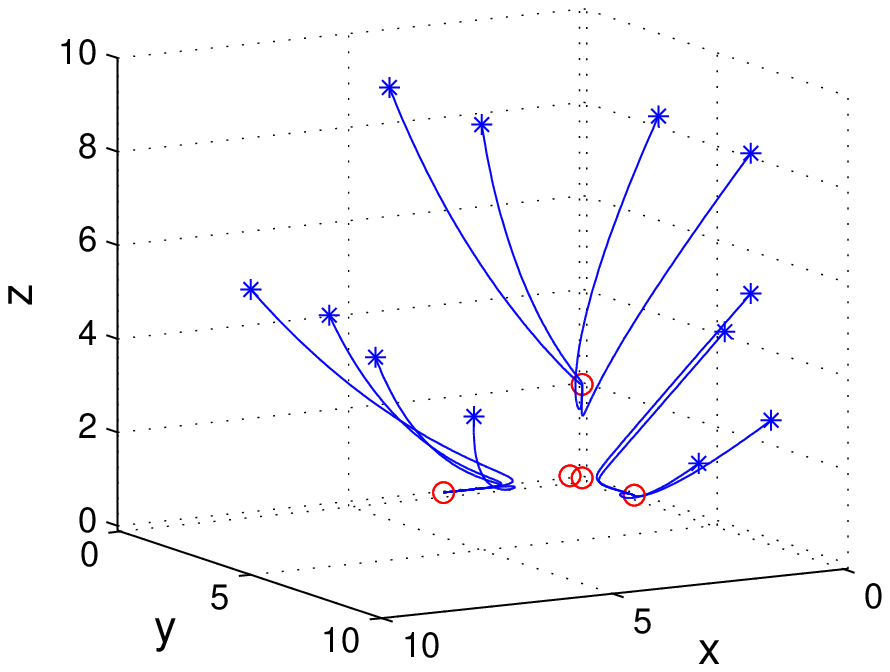}
  \caption{Example of trajectories for the model problem (\ref{model3d})
originating from different initial conditions and
converging to two different equilibria (top) and three equilibria (bottom).}
\label{ritratto_fase_3D}
\end{center}
\end{figure}

In Figure \ref{ritratto_fase_3D} (top) we show trajectories starting from the initial conditions $\boldsymbol{x}_{1}=(7,8,4)$, $\boldsymbol{x}_{2}=(8,7,10)$, $\boldsymbol{x}_
{3}=(8,7,4)$, $\boldsymbol{x}_{4}=(7,8,10)$, $\boldsymbol{x}_{5}=(5,8,4)$, $\boldsymbol{x}_{6}=(6,7,10)$, $\boldsymbol{x}_{7}=(6,7,4)$ and $\boldsymbol{x}_{8}
=(5,8,10)$, and converging to the point $E_3$ of coordinates $\left(0, 0, 9.5 \right)$ and to the point $E_4$ of coordinates $\left(  0.8511, 0.1489, 0 \right)$.

In Figure \ref{ritratto_fase_3D} (bottom) we show trajectories starting from the initial conditions $\boldsymbol{x}_{1}=(2,10,6)$, $\boldsymbol{x}_{2}=(4,10,10)$, $\boldsymbol{x}_
{3}=(8,10,4)$, $\boldsymbol{x}_{4}=(2,8,2)$, $\boldsymbol{x}_{5}=(9,8,5)$, $\boldsymbol{x}_{6}=(2,10,9)$, $\boldsymbol{x}_{7}=(1,9,3)$, $\boldsymbol{x}_{8}=(10,8,6)$, $\boldsymbol{x}_{9}=(5,5,9)$, $\boldsymbol{x}_{10}=(7,5,10)$, $\boldsymbol{x}_{11}=(2,5,9)$ and $\boldsymbol{x}_{12}=(10,5,6)$ and converging to the point $E_1$ of coordinates $\left(3, 0, 0 \right)$, $E_2$ of coordinates $\left(0, 2, 0 \right)$ and to the point $E_3$ of coordinates $\left( 0,0, 2\right)$.

To determine the separatrix surfaces for (\ref{model3d}), we need to consider a set of points as initial conditions in
a cubic domain $[0,\gamma]^3$, where $\gamma \in \RR^+$ (in the following we will fix $\gamma=10$).
Then, as in the 2D case, we take points in pairs and we check if trajectories of the two points converge to different
equilibria. If this the case, we proceed with a bisection algorithm to determine a separatrix
point. The bisection algorithm now is different, since the system presents three stable equilibria. In fact, the trajectories can evolve toward three different equilibrium points. As a consequence, we have to distinguish the points lying on the surface separating the trajectories evolving toward $E_1$ from the trajectories that evolve toward $E_2$ or $E_3$, the surface separating the initial conditions evolving toward $E_2$ from the initial conditions that evolve toward $E_1$ or $E_3$ and the surface separating the trajectories evolving toward $E_3$ from the trajectories that evolve toward $E_1$ or $E_2$.

More precisely, for the detection of separatrix points in the 3D case we use a technique very similar to that used in the 2D case.
At first we construct a grid on the faces of the cube and we apply the bisection algorithm  with the following initial conditions: 
$$
(x_{i_1},y_{i_2},0), \quad  (x_{i_1},y_{i_2},\gamma), \quad i_1=1, \ldots, n,  i_2=1, \ldots, n,
$$
$$
(x_{i_1},0,z_{i_2}), \quad (x_{i_1},\gamma,z_{i_2}), \quad i_1=1, \ldots, n,  i_2=1, \ldots, n,
$$
$$
(0,y_{i_1},z_{i_2}), \quad (\gamma,y_{i_1},z_{i_2}), \quad i_1=1, \ldots, n,  i_2=1, \ldots, n.
$$
The $N$ points found by the bisection algorithm are organized in a matrix 
$A=(a_{j,k})$, $j=1, \ldots, N$, $k=1,2,3$. 

\begin{figure}[ht!]
\begin{center}
  \includegraphics[height=.27\textheight]{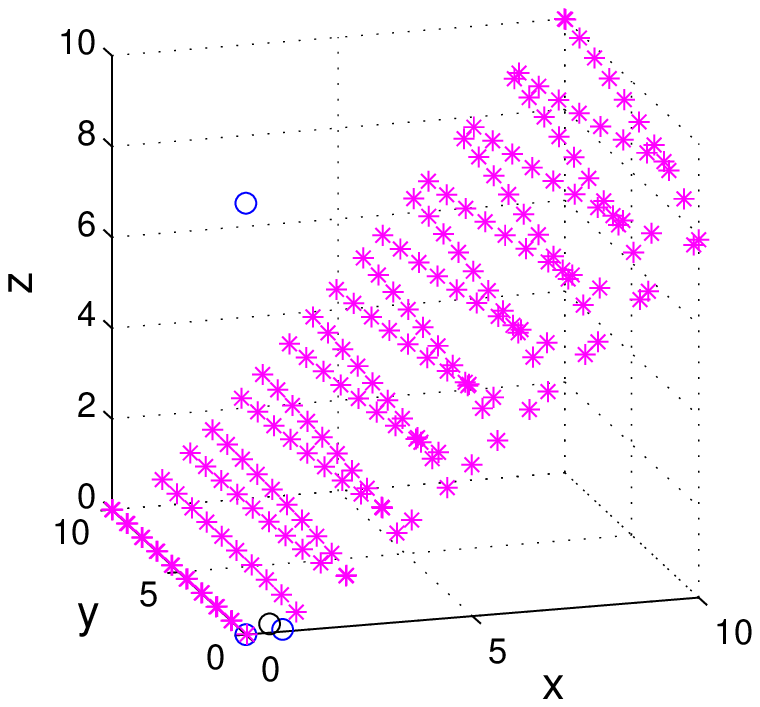}
  \includegraphics[height=.27\textheight]{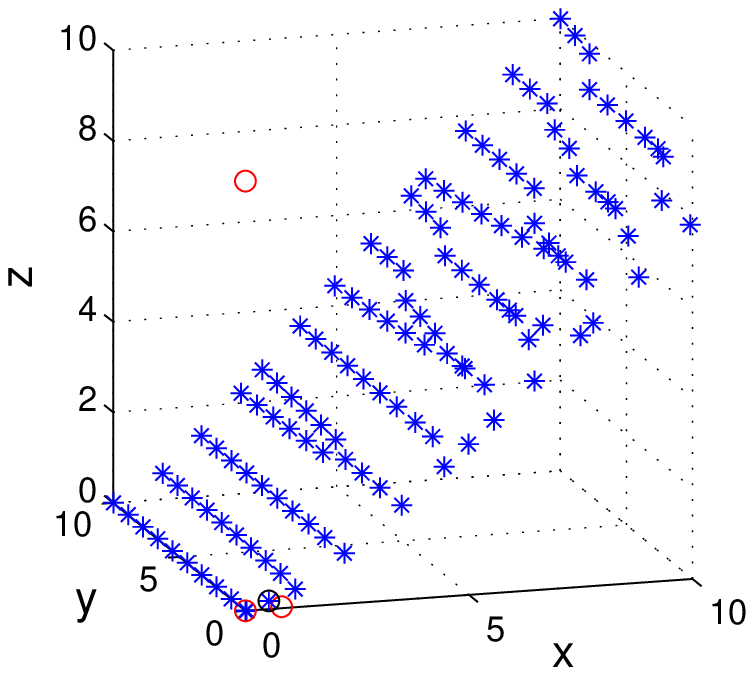} 
    \caption{Set of points detected by the bisection algorithm (top) and
  set of points found by the refinement algorithm (bottom) in the case of two stable equilibria.}
\label{figura3}
\end{center}
\end{figure}

\begin{figure}[ht!]
\begin{center}
  \includegraphics[height=.27\textheight]{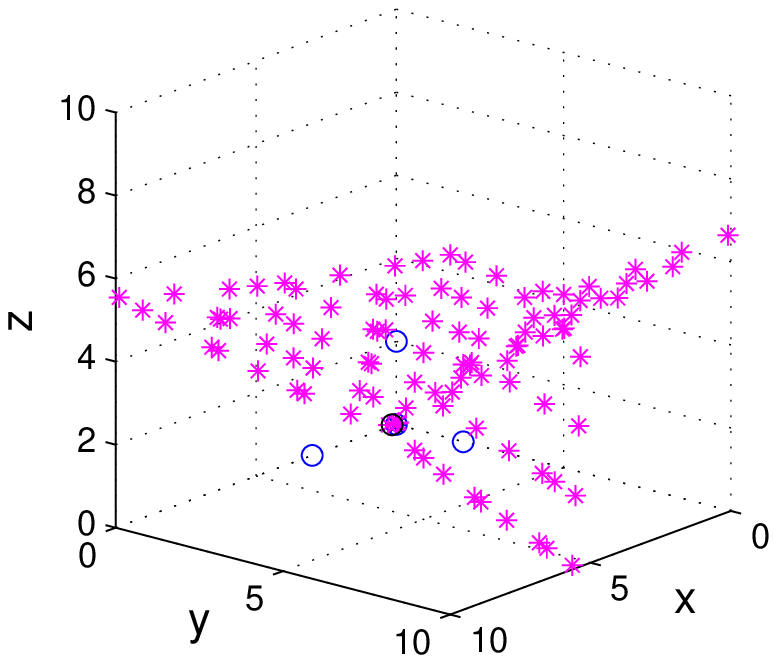}
  \includegraphics[height=.27\textheight]{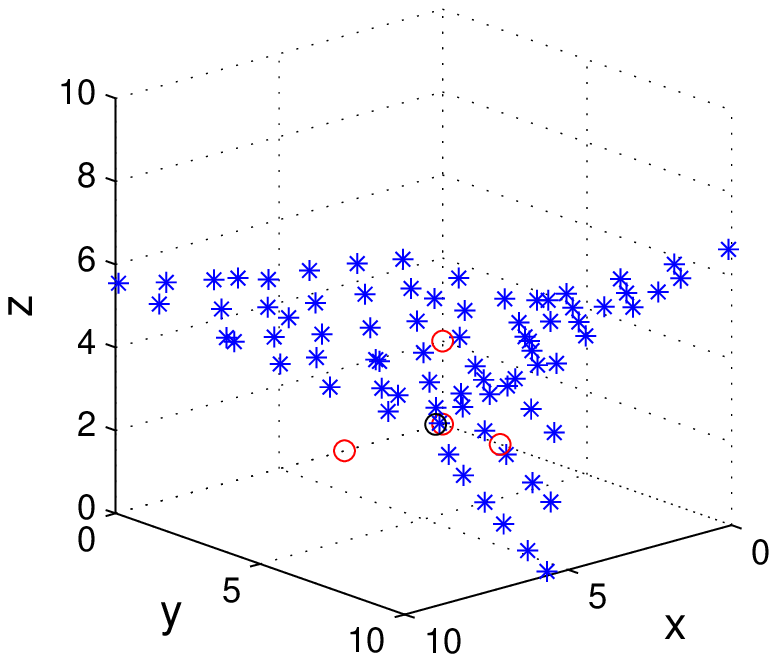} 
    \caption{Sets of points detected by the bisection algorithm (top) and
  sets of points found by the refinement algorithm (bottom) in the case of three stable equilibria.}
\label{figura4}
\end{center}
\end{figure}

As an example, in Figure  \ref{figura3} (top) and in Figure  \ref{figura4} (top) we show the points found by the bisection algorithm, in the case of two and three stable equilibria, choosing $n=10$ and $n=7$, respectively. From  Figure  \ref{figura4} (top), we deduce that we have to reconstruct separately the surface that determines the basin of attraction of $E_3$ and the surface that separates the trajectories tending to $E_1$ or $E_2$.
 For this aim we consider, starting from the matrix of points found by the bisection algorithm $A=(a_{j,k})$, $j=1, \ldots, N$, $k=1,2,3$, two submatrices $A^{'}=(a^{'}_{j,k})$, $j=1, \ldots, N^{'}$, $k=1,2,3$, and $A^{''}=(a^{''}_{j,k})$, $j=1, \ldots, N^{''}$, $k=1,2,3$. The points lying on the surface that determines the basin of attraction of $E_3$ are organized in the matrix $A^{'}$, while the remaining points that separate the trajectories tending to $E_1$ or $E_2$ are organized in the matrix $A^{''}$. 

To obtain smaller sets of nodes well distributed on the separatrix surfaces, we can proceed as follows. Let $B_{j,k}$, $j=1, \ldots, M$, $k=1,2,3$, a general matrix containing separatrix points lying on a separatrix surface.
We define
$$
M_x= \max_{ \displaystyle j} (b_{j,1}), \quad j=1, \ldots, M,
$$
$$
M_y= \max_{ \displaystyle j} (b_{j,2}), \quad j=1, \ldots, M,
$$
$$
M_z= \max_{ \displaystyle j} (b_{j,3}), \quad j=1, \ldots, M,
$$\\
and we divide the intervals $[0,M_x]$, $[0,M_y]$ and $[0,M_z]$ in $L$ subintervals.
We consider the following equispaced vectors in the intervals $[0,M_x]$, $[0,M_y]$ and $[0,M_z]$, respectively,
$x_l,$ $l=1, \ldots, L+1,$ $y_h,$ $ h=1, \ldots, L+1,$ $z_p,$ $p=1, \ldots, L+1,$
and we define
$$\\
I_{lhp}= \{j : b_{j,1} \in [x_l,x_{l+1}] 
$$
$$\\
b_{j,2} \in  [y_{h},y_{h+1}]  
$$
$$\\
  b_{j,3} \in  [z_{p},z_{p+1}] \},
$$
with $l=1, \ldots , L,$ $ h=1, \ldots, L$, $ p=1, \ldots, L$. Starting from the matrix $B=(b_{j,k})$ we find the matrix of the refined points $B^{'}=(b^{'}_{j,k})$, whose entries are given by: 
$$
b^{'}_{j,1}= \frac{ \sum_{j \in I_{lhp}} b_{j,1} }{ Card(I_{lhp})},  \quad l,h,p=1, \ldots , L,
$$
$$
b^{'}_{j,2}= \frac{ \sum_{j \in I_{lhp}} b_{j,2} }{ Card(I_{lhp})}, \quad l,h,p=1, \ldots , L, 
$$
$$
b^{'}_{j,3}= \frac{ \sum_{j \in I_{lhp}} b_{j,3} }{ Card(I_{lhp})}, \quad l,h,p=1, \ldots , L, 
$$
$j=1, \ldots,K$, where $K$ is the number of subintervals containing at least a point and $Card$ is the cardinality of the sets.

In the case of two equilibria we apply the refinement algorithm to the matrix $A$.
While, in the case of three equilibria, we obtain two different sets of points lying on the two different surfaces and we refine both sets, i.e. we apply the refinement algorithm to the matrices $A^{'}$ and $A^{''}$. 
These points will then be interpolated to reconstruct the required surfaces.

As an example, in Figure  \ref{figura3} (top) we show the points found by the bisection algorithm lying on the surface that separates the trajectories tending to $E_3$ or $E_4$ choosing $n=10$, in the case of two stable equilibria. The $N=195$ points have been refined taking  $L=13$. In this way, as shown in Figure \ref{figura3} (bottom), we obtain $K=127$ points.

In Figure  \ref{figura4} (top) we show the $N=102$  points found by the bisection algorithm, choosing $n=7$, in the case of three stable equilibrium points. $N^{'}=81$ points lie on the surface that determines the basin of attraction of the stable equilibrium point $E_3$ and $N^{''}=21$ points lie on the  surface separating  the trajectories tending to $E_1$ from those tending to $E_2$. In Figure \ref{figura4} (bottom) we show the points found
by the refinement algorithm taking $L=13$. We obtain $K^{'}=61$ points lying on the surface that determines the
trajectories tending to $E_3$ and $K^{''}=16$ points lying on the other surface. 

\section{Numerical experiments} \label{num_res}

\begin{figure}[ht!]
\begin{center}
  \includegraphics[height=.26\textheight]{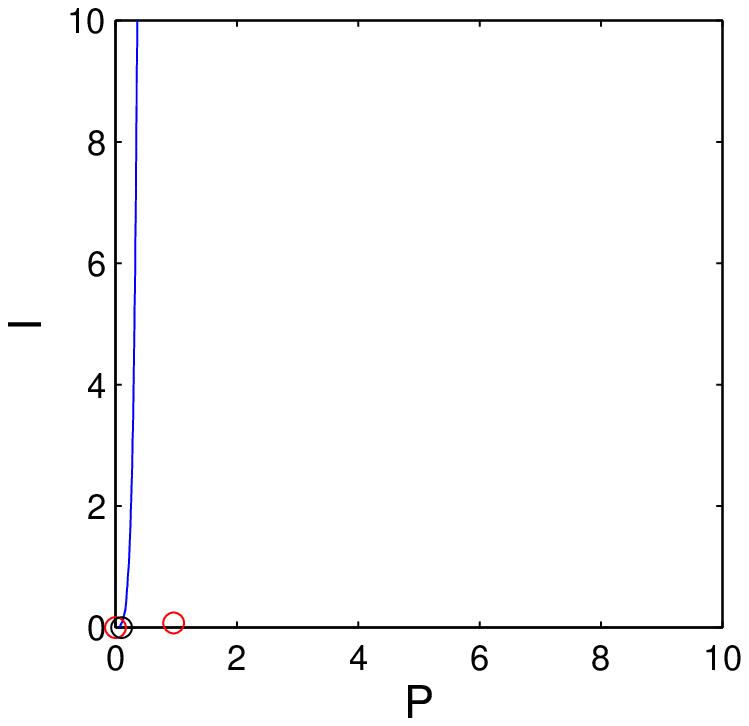} \\
  \includegraphics[height=.24\textheight]{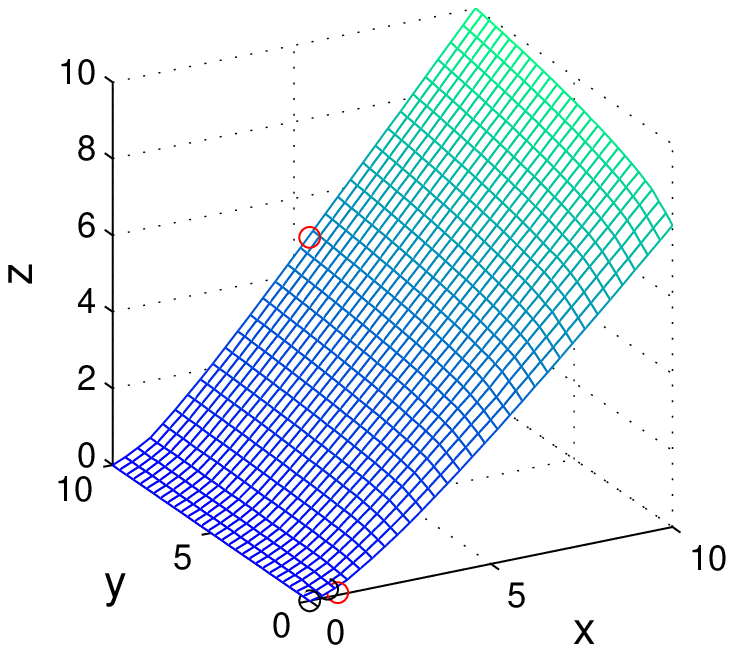} 
  \caption{Approximation of separatrix curve using the Gneiting's $C^2$ (top) and approximation of separatrix surface  using the Wendland's $C^2$  (bottom).}
\label{figura6}
\end{center}
\end{figure}

\begin{figure}[ht!]
\begin{center}
  \includegraphics[height=.24\textheight]{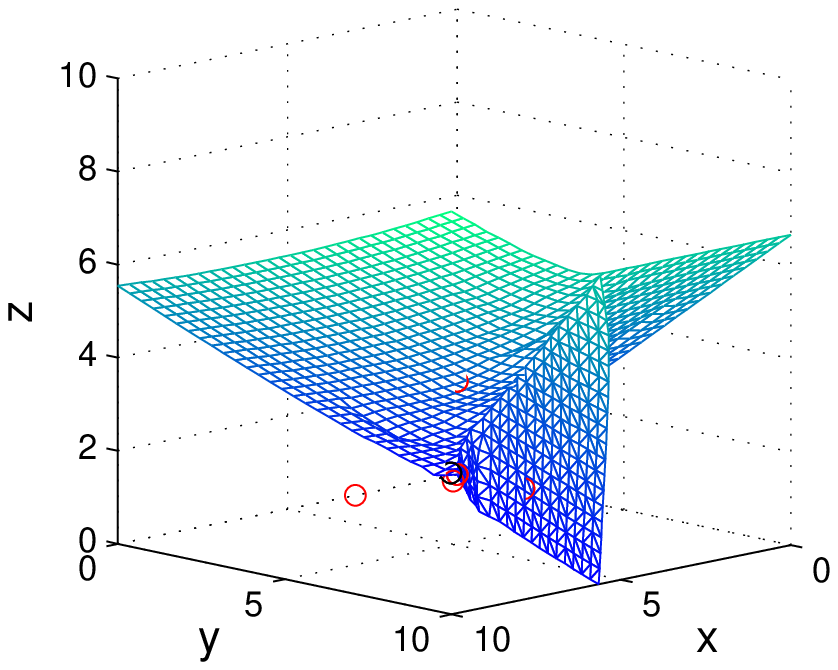} \\
  \includegraphics[height=.24\textheight]{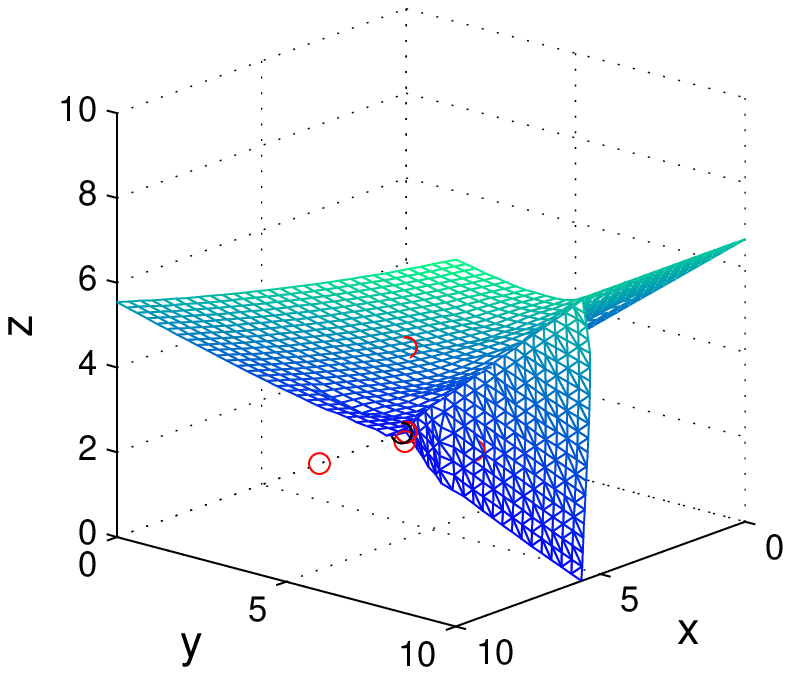} 
  \caption{Approximation of separatrix surfaces using the Wendland's $C^2$ function (top) and the Wu's $C^4$  function (bottom).}
\label{figura7}
\end{center}
\end{figure}

In this section we summarize the extensive experiments performed to test our detection and approximation algorithms.
As far as the accuracy of the Partition of Unity method,
a crucial task concerns the choice of the shape parameter $c$ of the compactly supported function. In fact, it can significantly affect the approximation result and, therefore, the quality of the separatrix curves and surfaces.
From our study we found,  that good shape parameter values for all the functions listed in Section 2 are in the ranges $0.01 \le c \le 0.05 $ (case 2D) and $0.001 \le c \le 0.01$ (case 3D). In Figure \ref{figura6} (top)  we show the curve obtained approximating the refined data set when we consider the value $c=0.015$ as shape parameter for the  Gneiting's $C^2$ function $\tau_1$  and a number $d=4$ of partitions of $\Omega$.  Figure \ref{figura6} (bottom) shows bistability in the 3D case: the separating surface is reconstructed using  $c=0.005$ for the  Wendland's $C^2$ function $\varphi_1$  and $d= 4$  partitions of $\Omega$.
In Figure \ref{figura7}  the two surfaces partitioning the domain in three regions are shown. For the surface that determines the basin of attraction of $E_3$ we consider, for both the Wendland's $C^2$ and the Wu's $C^4$ functions,   the value $c=0.005$  and a number $d= 4$ of partitions of $\Omega$.
 In order to reconstruct the surface that separates the paths tending to $E_1$ or $E_2$, we interpolate  the  points found by the bisection algorithm exchanging the $x$ axis with the $z$ axis. Acting in this way we approximate the surface on a triangular domain in the $xz$ plane. We consider, for both the Wendland's $C^2$ and the Wu's $C^4$ functions, the value $c=0.005$ as shape parameter  and a number $d=3$ of partitions of $\Omega$.

	Our results, shown in this section turn out to be accurate and moreover also  stable, since compactly supported RBFs are used. It is well known that the use of the latter leads to a good compromise between accuracy and stability.

\section{Summary of results and future work} \label{concl}

In this paper we presented an approximation method for the detection of points lying on
the separatrix curve for the model (\ref{model2d_hilker}) and the separatrix surfaces for the model (\ref{model3d}).

An efficient algorithm based on the Partition of Unity method, which uses  complactly supported radial basis  functions as
local approximants, is used for the reconstruction of separatrix curves and surfaces.
It was already used in previous papers (see \cite{cavoretto11,C-D-P-V}), but here we considered
an extension and a refinement to account for a different model. In particular, 
the approximation scheme has been improved as far as portability is concerned, in that now it
works also for a dynamical system of dimension three with three stable equilibria.

Even if the algorithm for the detection of separatrix points in case of three equilibria is performed without any specific intervention of the user on the code, the interpolation instead needs a direct manipulation of the points to be interpolated, as stressed in Section \ref{num_res}. Thus work in progress consists in further investigations about the numerical interpolation scheme in order to 
reconstruct separatrix surfaces defined in arbitrary domains or also implicitly defined without any treatment by the user.


\section*{Acknowledgements}
This research was partially supported by the project ``Metodi numerici nelle scienze applicate'' of the Dipartimento di Matematica ``Giuseppe Peano''.

\end{document}